\documentclass[12pt, reqno]{amsart}
\usepackage{ifthen,amsfonts,amsmath,amssymb,epic,eepic}
\usepackage{epsfig,color}

\makeatletter
\@addtoreset{equation}{section}
\makeatother

\theoremstyle{definition}

\def\fnum{equation}
\newtheorem{Thm}[\fnum]{Theorem}
\newtheorem{Cor}[\fnum]{Corollary}

\newtheorem{Lem}[\fnum]{Lemma}

\newtheorem{Def}[\fnum]{Definition}

\newtheorem{Rem}[\fnum]{Remark}
\newtheorem{Pro}[\fnum]{Proposition}

\numberwithin{equation}{section}

\newcommand{\eps}{{\varepsilon}}

\newcommand{\nn}{{\bf{n}}}
\newcommand{\Ric}{{\text{Ric}}}
\newcommand{\Scal}{{\text{Scal}}}

\def\RR{{\bf  R}}
\def\ZZ{{\bf  Z}}
\def\SS{{\bf  S}}

\newcommand{\cG}{{\mathcal{G}}}
\newcommand{\cF}{{\mathcal{F}}}

\newcommand{\cL}{{\mathcal{L}}}

\newcommand{\eqr}[1]{(\ref{#1})}

\begin{document}

\title[Singular limit laminations, Morse index, and positive scalar curvature]
{Singular limit laminations, Morse index, and positive scalar curvature}

\author{Tobias H. Colding}%
\address{Courant Institute of Mathematical Sciences and Princeton University\\
251 Mercer Street\\ New York, NY 10012 and Fine Hall, Washington
Rd., Princeton, NJ 08544-1000}
\author{Camillo De Lellis}%
\address{Max-Planck-Institute for Mathematics in the Sciences\\
Inselstr. 22 - 26, 04103 Leipzig / Germany
}

\thanks{The first author was partially supported by NSF Grant DMS 0104453.}

\email{colding@cims.nyu.edu and delellis@mis.mpg.de}

\begin{abstract}
For any $3$-manifold $M^3$ and any nonnegative integer ${\bf{g}}$, 
we give here 
examples of metrics on $M$ each of 
which has a sequence of 
embedded minimal surfaces of genus ${\bf{g}}$ and without Morse 
index bounds (all our surfaces will be orientable).  
On any spherical space form $\SS^3/\Gamma$  we construct such a 
metric with positive scalar curvature.
More generally we construct such a metric with $\Scal >0$ 
(and such surfaces) on any $3$-manifold which carries a 
metric with $\Scal>0$.
\end{abstract}

\maketitle

\section{Introduction}

For any $3$-manifold $M^3$ and any nonnegative integer ${\bf{g}}$, 
we give here 
examples of metrics on $M$ each of 
which has a sequence of 
embedded minimal surfaces of genus ${\bf{g}}$ and without Morse 
index bounds (all our surfaces will be orientable).  
On any spherical space form $\SS^3/\Gamma$  we construct such a 
metric with positive scalar curvature.
More generally we construct such a metric with $\Scal >0$ 
(and such surfaces) on any $3$-manifold which carries a 
metric with $\Scal>0$; see Theorem 
\ref{t:gensc} below.
In all but one of our examples the Hausdorff limit will be a 
singular minimal lamination.  The singularities being in each 
case exactly two 
points lying on a closed leaf (the leaf is a strictly stable sphere). 
 
There are two prior examples of embedded minimal surfaces in $3$-manifolds 
without Morse index bounds.  In \cite{CH1} it was shown that even in one 
dimension less (i.e., for simple closed geodesics on surfaces) there are 
examples of metrics without Morse index bounds.  \cite{CH1} also gave 
examples on any $3$-manifold of a metric which has 
embedded minimal tori without such bounds.  
In \cite{HaNoRu} examples were given of metrics on any $M^3$ that have 
embedded minimal spheres without bounds.  As mentioned above 
in this paper we are not only interested in giving such examples for any 
genus and of metrics with positive scalar curvature but also 
in a particular type of degeneration of the surfaces.  

We use in part ideas 
of Hass-Norbury-Rubinstein 
\cite{HaNoRu} to achieve this (and in the process answer a 
question of theirs).  As in \cite{HaNoRu}, but  
unlike the examples in \cite{CH1}, 
the surfaces will have no uniform curvature bounds. 
In fact, it follows easily (see appendix  B of \cite{CM4}) that 
 if
$\Sigma_i\subset M^3$ is a sequence of embedded  minimal surfaces with 
uniformly
bounded curvatures, then a subsequence  converges to a smooth lamination.  
Moreover, with the right notion of being generic,  
the following seems likely (by \cite{CH1} bumpy is 
not the right generic notion):

\vskip2mm
\noindent
{\bf{Conjecture}}:  Let $M^3$ be a closed $3$-manifold with a 
generic metric and $\Sigma_i\subset M$ a sequence of embedded minimal 
surfaces of a given genus.  If any limit 
of the $\Sigma_i$'s is a 
{\underline{smooth}} (minimal) lamination, then the sequence $\Sigma_i$ has a  
uniform Morse index bound.
\vskip2mm

A codimension one {\it lamination} of $M^3$ is a collection $\cL$ of
smooth disjoint connected surfaces (called leaves)
such that  $\cup_{\Lambda \in \cL} \Lambda$ is closed.
Moreover, for each $x\in M$ there exists an
open neighborhood $U$ of $x$ and a local coordinate chart, $(U,\Phi)$, with
$\Phi (U)\subset \RR^3$
such that in these coordinates the leaves in $\cL$
pass through the chart in slices of the
form $(\RR^2\times \{ t\})\cap \Phi(U)$.

A lamination is said to be minimal if the leaves are (smooth)
minimal surfaces.  If the union of the leaves  is all of $M$, then it
is a foliation.  

There are two results that support this conjecture.  
The first concerns the corresponding conjecture 
in one dimension less (that is for geodesics on surfaces); see
\cite{CH2}, \cite{CH3}. 
The second concerns the conjecture for 
$3$-manifolds with 
positive scalar curvature; see \cite{CM3}.  However, there are
examples where the limit is not smooth  
as the following shows:  

\begin{Thm}   \label{t:genus}
On any $3$-manifold, $M^3$, and for any nonnegative integer $\bf{g}$, 
there exists a metric and a sequence of embedded minimal  
surfaces of 
genus $\bf{g}$ with Morse index going to infinity and converging to a 
singular (minimal) lamination $\cL$.  This can be done so that the singular 
set of $\cL$ consists of two points lying on a leaf which is a strictly 
stable $2$-sphere.
\end{Thm}

\vskip2mm
For manifolds which carry a metric with positive scalar curvature we 
use a connected sum construction to show (cf. section 5 of 
Gromov-Lawson \cite{GrLa} and theorem 4 of Schoen-Yau \cite{ScYa}):

\begin{Thm}   \label{t:gensc}
(See fig. \ref{fig8}).  
Any $3$-manifold which carries a metric with positive scalar curvature 
has for any nonnegative integer ${\bf{g}}$ 
a metric with positive scalar curvature and a sequence of embedded 
minimal surfaces of genus ${\bf{g}}$ as in Theorem \ref{t:genus}.
\end{Thm}

As a consequence we get by \cite{GrLa}, \cite{ScYa}:

\begin{Cor}   \label{c:gensc}
Any manifold of the form
\begin{equation}
\SS^3/\Gamma_1\#\cdots\#\SS^3/\Gamma_k\#
\SS^2\times \SS^1\#\cdots\#\SS^2\times \SS^1\, ,
\end{equation}
where $\SS^3/\Gamma_i$ is a spherical space form, has for any nonnegative 
integer $\bf{g}$ a metric with positive scalar curvature and a sequence 
of embedded 
minimal surfaces of genus ${\bf{g}}$ as in Theorem \ref{t:genus}.
\end{Cor}

\vskip2mm
The following is a different kind of example (different 
from \cite{CH1}; however not bumpy) 
that illustrates why generic is needed in the above conjecture:  

\begin{Thm}  \label{t:product}
In $\SS^2\times \SS^1$ with the product metric, there 
is a sequence of embedded minimal tori with Morse index going to 
infinity.  Moreover, these converge to the foliation by 
parallel $\SS^2\times \{t\}$.  
\end{Thm}

The next four sections contain the proofs of the above three
theorems.  In Section \ref{s:metrics} we show how to generalize Theorems
\ref{t:genus} and \ref{c:gensc} to where the singular set contains points 
on any given
finite collection of disjoint embedded strictly stable $2$-spheres.  Finally,
in Section \ref{s:speculation} we return to a result shown in Section 
\ref{s:neck}
and speculate on how the space of noncompact embedded minimal annuli
limiting a strictly stable $2$-sphere look like.   Moreover,
we speculate there on what the
structure of this space of annuli might imply for structure of 
the singular set of a
limit lamination for a generic metric.

\vskip2mm
Recall that if $\Sigma^2\subset M$
is a closed minimal surface, then the {\it Morse index} of
$\Sigma$ is the index of the
critical point $\Sigma$ for the area functional, i.e., the number of
negative
eigenvalues (counted with multiplicity) of the second
derivative of area.  If $\Sigma$ has a unit normal $\nn$,
 the second derivative of area at $\Sigma$
in the direction of a normal variation $u\,\nn$ is
$-\int_{\Sigma}u\,L\,u$ where
$L \,u= \Delta\,u + [|A|^2+\Ric_M(\nn,\nn)]\,u$;
so the Morse index is the number of
negative eigenvalues of $L$.
(By convention,  an
eigenfunction $\phi$ with eigenvalue
$\lambda$ of $L$ is a solution of
$L\,\phi+\lambda\, \phi=0$.)
$\Sigma$ is said to be
{\it stable} if the index is zero.  
A metric on $M^3$ is {\it bumpy} if
each closed minimal surface is a nondegenerate critical point,
i.e., $L \,u = 0$ implies $u\equiv 0$. By a 
result of B. White bumpy
metrics are generic; that is the set of bumpy metrics
contain a countable intersection of open dense subsets.  We use
throughout the normalization of the curvature so that the round unit
$3$--sphere has sectional curvature $1$ and scalar curvature $3$.  

Our interest in whether the Morse index is bounded 
for embedded minimal tori in a $3$-manifold comes in part 
from its connection with the spherical space form problem; 
see \cite{PiRu}, \cite{CM2}.

\vskip2mm
Part of this work was done while the first author was visiting Scuola Normale 
Superiore in Pisa.  He wishes to thank Mariano Giaquinta and Carlo 
Mantegazza for making the visit possible.

\section{The metric and surfaces near the stable $2$-sphere}
\label{s:neck}

Following \cite{HaNoRu} (see also \cite{HsLa}) we look at metrics on 
$\SS^2\times \RR$ of the form 
\begin{equation}   \label{e:metric}
ds_0^2=dr^2+\lambda^2(r)\,(d\phi^2+\sin^2 \phi\,d\theta^2)\, .
\end{equation}
Here $(\phi,\theta)$ are spherical coordinates on $\SS^2$ and $r\in \RR$.    
Computing the scalar curvature of 
the warped product gives
\begin{equation}  \label{e:scalwrap}
\Scal_M=-2\frac{\lambda''}{\lambda}+\frac{1-(\lambda')^2}{\lambda^2}\, .
\end{equation}
To find our minimal surfaces we 
consider on the infinite strip $[0,\pi]\times\RR$ the degenerate metric
\begin{equation}   \label{e:cube}
ds^2=\lambda^2(r)\,\sin^2 \phi \,(dr^2+\lambda^2(r)\,d\phi^2)
\end{equation} 
and calculate the geodesics in this metric.  Our minimal surfaces will be 
the preimages of simple closed geodesics in the metric \eqr{e:cube} under 
the map $(\phi,\theta,r)\to (\phi,r)$.  For completeness we will now see why 
these preimages are minimal.  So let $\Sigma$ be a surface of 
the form $\SS^1\times \gamma$ where $\gamma (t)= (\phi (t), r(t))$ 
is a curve in $[0,\pi]\times \RR$ (below $\phi$ will be different from 
$0$ and $\pi$ for the curve $\gamma$ so the preimage of each 
$\gamma (t)$ is indeed a circle). A surface is minimal if and only if 
the first variation of the area functional is zero for any smooth 
vector field perpendicular to it. Since the rotations 
$\theta\,\to\, \theta\, +$ constant preserve the metric \eqr{e:metric} 
and $\Sigma$, 
it is sufficient to check that the first variation vanishes with 
respect to vector fields invariant for this family 
of isometries. Being perpendicular to $\Sigma$, these 
vector fields are of the form 
$v=v_\phi (\phi, r)\, \partial_\phi + v_r (\phi, r) \,\partial_r$.  
Thus checking first variation of the area for $\Sigma$ 
is equivalent to check the first variation of the functional
\begin{equation}
  F (\gamma)=\int_{\gamma}\, \textrm{length} ( \SS^1 \times \{\gamma
  (t)\})=\int_{\gamma}\, 2\pi\, \lambda (t)\, \sin (\phi (t))
\end{equation}
in the space of curves of $[0,\pi]\times \RR$ with the metric 
$dr^2+\lambda^2 (r)\, d\phi^2$. Notice that $F(\gamma)$ is $2\pi$
times the length of 
$\gamma$ in the metric \eqr{e:cube} and hence the first variation of $F$ 
vanishes if and only if $\gamma$ is a geodesic in \eqr{e:cube}.

For a unit speed geodesic in \eqr{e:cube} (throughout this paper all
geodesics will have unit speed) 
\begin{align}  
r''&=-2\,\frac{\cos \phi}{\sin \phi}\,r'\,\phi'
-\frac{\lambda'(r)}{\lambda (r)}\,(r')^2
+2\,\lambda'(r)\,\lambda (r)\,(\phi')^2\, ,\label{e:r''}\\
(r')^2+\lambda^2\,(\phi')^2&=\lambda^{-2}\,\sin^{-2}\phi\, .\label{e:speed}
\end{align}
From \eqr{e:r''} it follows that 
if $\lambda'\geq 0$, then provided $r'>0$
\begin{equation}  \label{e:concl}
\frac{d}{dt}\log r'\geq \frac{d}{dt}\log \sin^{-2} \phi
-\frac{\lambda' (r)}{\lambda (r)}\,r'\, .
\end{equation}
In particular \eqr{e:concl} yields that if $r'(0)>0$, then $r'(t)>0$ 
for all $t>0$.   Namely, suppose that $r'(t_0)=0$ and that 
$t_0=\inf \{t>0\,|\,r'(t)=0\}$, applying \eqr{e:concl} yields a 
contradiction.  It follows that if $r'(0)>0$, 
then the geodesic is simple.   Moreover, integrating \eqr{e:concl} yields 
for $t_2>t_1$ 
\begin{equation}  \label{e:repelling}
\frac{r'(t_2)}{r'(t_1)}
\geq \frac{\sin^{2}\phi (t_1)}{\sin^{2}\phi (t_2)}\,
\exp \left( C_1\,(r(t_1)-r(t_2))+C_2\right)\, .
\end{equation}   
One may also easily check that if $\lambda''> 0$, 
then the only curve where $r$ is constant that is a geodesic is for $\{r=0\}$ 
(this follows for instance since the only level set of $r$ in \eqr{e:metric}  
that is a minimal surface is $\{r=0\}$).  
Finally, it follows from \eqr{e:speed} and \eqr{e:repelling} 
that the boundary of 
the infinite strip is repelling.  
That is, the only geodesics that intersect the
boundary are $\{r=0\}$, $\{\phi=\pi\}$, and $\{\phi=0\}$.
Using this we can now show:  

\begin{Pro}  \label{p:neck}
For $\eps>0$ set $\lambda_{\eps}(r)=\cosh (\eps\,r)$.  
On $M=\SS^2\times_{\lambda_{\eps}} \RR$ (see fig. \ref{fig1}), $\SS^2\times \{0\}$ 
is the only closed minimal surface and there is a singular minimal 
lamination $\cL$ on $M$ with antipodal points 
on $\SS^2\times\{0\}$ as the only singularities of $\cL$; see fig. \ref{fig2}.  
Moreover, there 
is a sequence of embedded minimal annuli $\Sigma_i$ with $\Sigma_i\to \cL$ 
and with Morse index going to 
infinity.
\end{Pro}

\begin{figure}[htbp]
\begin{center}
    \begin{minipage}[t]{0.45\linewidth}
    \centerline{\input{pic1.pstex_t}}
    \caption{The strictly stable $2$--sphere in the warped metric
    \eqref{e:metric} with $\lambda (r)=\cosh (\eps r)$}
    \label{fig1}
    \end{minipage}
    \begin{minipage}[t]{0.45\linewidth}
    \centerline{\input{pic2.pstex_t}}
    \caption{The singular lamination in half of a neighborhood of the strictly
    stable $2$--sphere}
    \label{fig2}
    \end{minipage}
\end{center}
\end{figure}

\begin{proof}
To prove this proposition all we need is to find the corresponding 
geodesic lamination and simple geodesics 
on the infinite strip $[0,\pi]\times\RR$ with the 
degenerate metric (see fig. \ref{fig3})
\begin{equation}   \label{e:cube2}
ds^2
=\cosh^2(\eps\,r)\,\sin^2 \phi \,(dr^2+\cosh^2(\eps\,r)\,d\phi^2)\, .
\end{equation} 
Let $\gamma_{\delta}(t)=(\phi_{\delta}(t),r_{\delta}(t))$ 
be a geodesic in the metric \eqr{e:cube2} with 
$(\phi_{\delta}(0), r_{\delta}(0))=(\pi/2,0)$ and so 
that the angle between $\gamma_{\delta}'(0)$ 
and $\{r=0\}$ 
is $\delta$.  We extract a sequence $\gamma_{i}=\gamma_{\delta_i}$ 
with $\delta_i\to 0$ and 
which converges in the Hausdorff sense.
By the equations for geodesics above it follows 
that every $\gamma_{i}$ is simple and that 
$\gamma_{i}\to \cG$ as 
$i\to \infty$, where $\cG$ is a geodesic lamination consisting of $\{r=0\}$ 
and two infinite 
geodesics $\gamma_\infty$ and $\gamma_{-\infty}$
which lie on each side of $\{r=0\}$ and spiral into it.  The surfaces 
$\Sigma_{i}$ and the singular minimal lamination $\cL$ 
can now be taken to be the preimages of the geodesics 
$\gamma_{i}$ and of $\cG$. 

It follows (by a standard argument) that for any $r_0>0$ 
the Morse index of $T_{r_0}(\SS^2\times \{0\})\cap \Sigma_{\delta}$ goes to 
infinity as $\delta\to 0$ (basically it follows easily, 
at least for $r$ small, that the preimage 
of each ``turn'' in $\gamma_{\delta}$, see fig. \ref{fig2} and fig. \ref{fig4}, 
corresponds to a small 
neck that contributes to the index). Alternatively we can use the fact
that Jacobi fields on the geodesics in \eqref{e:cube} lift to Jacobi
fields on the respective minimal surfaces in \eqref{e:scalwrap} and
then reason as in \cite{HaNoRu}.
Finally, it follows easily from the 
maximum principle, as in the proof of proposition 1.8 of \cite{CM3} 
(the sublevel sets $\{r\leq r_0\}$ are strictly mean convex for $r_0>0$), 
that $\SS^2\times \{0\}$ is the only closed minimal surface in $M$.  
\end{proof}

\begin{figure}[htbp]
\begin{center}
    \begin{minipage}[t]{0.45\linewidth}
    \centerline{\input{pic3.pstex_t}}
    \caption{The upper half--strip with the degenerate metric \eqref{e:cube} where $\lambda (r)=\cosh (\eps r)$}
    \label{fig3}
    \end{minipage}
    \hspace{1cm}
    \begin{minipage}[t]{0.45\linewidth}
    \centerline{\input{pic4.pstex_t}}
    \caption{The geodesic $\gamma_\delta$ in the same half--strip}
    \label{fig4}
    \end{minipage}
\end{center}
\end{figure}

We will later need to deal with that the geodesics 
$\gamma_{i}$ and $\gamma_\infty$
(and hence also the corresponding minimal surfaces)
cross in many points. In order to prove our theorems we will use that,
by the next lemma, we can 
choose the $\Sigma_i$'s and $\Sigma_{\infty}$ so 
that in a neighborhood of some point of
$\Sigma_\infty$ the $\Sigma_i$'s can be completed to a smooth
minimal foliation. 

\begin{Lem}\label{l:spirals}
Consider on $[0,\pi]\times \RR$ the degenerate 
metric \eqref{e:cube} where
$\lambda\in C^{1,1}$, $\lambda (r)=\lambda (-r)$, 
$\lambda'\geq 0$ on $[0, +\infty[$, 
and $\lambda (r)=\cosh (\eps r)$ for some $\eps>0$ 
in a neighborhood of $0$. 
For any fixed $\rho>0$, we can assume that
the geodesics $\gamma_i$, $\gamma_\infty$ 
constructed in Proposition \ref{p:neck} pass through $(\pi/2,
\rho)$. Thus $\{\gamma_i\}\cup\{\gamma^+\}$ can be
completed to a smooth geodesic foliation in a punctured ball centered at
$(\pi/2, \rho)$. 
\end{Lem}

\begin{proof} For any given $\gamma$
which starts at $\{\pi/2, 0\}$ and any integer $N$ let
\begin{eqnarray}
&&\alpha (\gamma) \mbox{ be the angle between $\{r=0\}$ and 
$\gamma$\, .}\nonumber\\
&&r_N (\gamma) \mbox{ be the $N$--th crossing between $\{\phi=\pi/2,
  r>0\}$ and $\gamma$}\, . \nonumber
\end{eqnarray}
Fix a $\gamma_\infty$ and a sequence $\gamma_i\to \gamma_\infty$ 
given by Proposition \ref{p:neck}. It is not difficult to see that for
sufficiently large $N$'s there exist $i$, $j$ such that 
$r_N (\gamma_j)\leq \rho\leq r_N (\gamma_i)$. Since
$\{\phi=\pi/2\}$ is a geodesic, any other geodesic
$\gamma$ crosses $\{\phi=\pi/2\}$ transversally. This easily implies 
that $r_N (\gamma)$ is a continuous function of the starting angle
$\alpha (\gamma)$.

Thus, varying this angle between $\alpha (\gamma_j)$ and $\alpha (\gamma_i)$,
we find a geodesic $\tilde{\gamma}_N$ starting at $(\pi/2, 0)$ with
$r_N (\tilde{\gamma}_N)=\rho$. 
Clearly $\alpha (\tilde{\gamma}_N)\to 0$ as $N\to \infty$. Hence,
we can extract a subsequence of $\{\tilde{\gamma}_N\}$ converging to
a geodesic lamination as in Proposition \ref{p:neck}. 
\end{proof}

The next definition and proposition are needed only in the proof of 
Theorem \ref{t:gensc}. 

\begin{Def}\label{d:parallel}
Let $z\in\Omega\subset \SS^3$ be an open subset of the round unit
$3$-sphere and suppose that $\cF$ is a foliation by great spheres 
of $\Omega$.  We say that the foliation is  {\em parallel} at 
$z$ if $\sup_{y\in \Lambda}
\text{dist}\, (y, \Lambda')=\text{dist}\, (z, \Lambda')$
where $\Lambda$, $\Lambda'\in \cF$ and $z\in \Lambda$ 
($\Lambda$ is said to be the {\it central} leaf of $\cF$). 
\end{Def}

This particular kind of foliation is needed in 
the proof of Theorem \ref{t:gensc} to make the connected sum construction. 

\begin{Pro}   \label{p:basiclam}
On $\SS^3$, 
there is a metric with $\Scal>0$ which has a singular 
lamination and a sequence of embedded 
minimal surfaces of genus $0$ as in Theorem \ref{t:genus}.
We can choose the metric so that these minimal
spheres can be completed in a neighborhood of a point $x$ to a 
foliation by great spheres parallel at $x$. Moreover, in an open
(nonempty) set disjoint from the minimal spheres 
the sectional curvature of the metric on $\SS^3$ is constant 1.
\end{Pro}

\begin{proof} Fix on $\SS^2\times \RR$ a metric with positive scalar
  curvature of the form
  \eqref{e:metric} where  $\lambda$ satisfies, 
  for some positive constants $a,b,c,\delta,\eps$,
  \begin{equation}
  \left\{\begin{array}{ll}
  \lambda (r)=\lambda (-r)\, ,& \\
  \lambda (r)=c \cosh (\eps r) & \mbox{in a neighborhood of $0$}\, ,\\ 
  \lambda' (r)\geq 0 & \mbox{for $r\in [0, \infty[$}\, ,\\ 
  \lambda(r)=1 & \mbox{for $r\in [a, \infty[$}\, ,\\
  \lambda (r)=\sin (r+\pi/2-a) & \mbox{for $r\in \, ]a-\delta, a[$}\, .
  \end{array}\right.
  \end{equation}
  $\lambda$ can be chosen $C^{1,1}$
  and $C^\infty$ on $\RR\setminus \{a,-a\}$. In particular,
  endowing $[0,\pi]\times \RR$ with the degenerate metric
  \eqref{e:cube}, by Lemma \ref{l:spirals} 
  there are geodesics $\gamma^+$, $\gamma^-$ through
  $(\pi/2, a)$ and $(\pi/2, -a)$, respectively, and spiraling into
  $\{r=0\}$. Moreover, again by Lemma \ref{l:spirals}, there is a
  sequence of geodesics $\gamma_i$ passing through $(\pi/2, a)$ and
  $(\pi/2, -a)$ which converges to the lamination 
  $\gamma^+\cup \gamma^-\cup \{r=0\}$.
  Define $\tilde{\lambda}$ by 
  \begin{equation}
  \tilde{\lambda} (r)=\left\{
  \begin{array}{ll}
  \lambda (r) & \mbox{for $r\in [-a,a]$}\, ,\\
  \sin (r+\pi/2-a) & \mbox{for $r\in [a, a+\pi/2]$}\, ,\\
  \sin (r+\pi/2+a) & \mbox{for $r\in [-a-\pi/2, -a]$}\, .
  \end{array}\right.
  \end{equation}
  Clearly $\tilde{\lambda}\in C^\infty$.
  On $\SS^2\times [-a-\pi/2, a+\pi/2]$ identify each of the spheres 
  $\SS^2\times \{a+\pi/2\}$ and $\SS^2
  \times \{-a-\pi/2\}$ to a point to get the smooth metric 
  \begin{equation}\label{e:finalmetric}
  dr^2 + \tilde{\lambda}^2 (r)\, (d\phi^2 + \sin^2 \phi\, d\theta^2)\,
  \end{equation}
  in $\SS^3$. This
  $\SS^3$ is obtained (loosely speaking) by capping off a neck with two standard
  half--$\SS^3$'s, $S^+$ and $S^-$.

  On $[0,\pi]\times [-a-\pi/2, a+\pi/2]$ with the degenerate metric
  \begin{equation}
  \tilde{\lambda}^2 (r) \sin^2 \phi \, (dr^2 + \tilde{\lambda}^2 (r)
  \, d\phi^2)\,
  \end{equation}
  the curve $\gamma^+\cap [0,\pi]\times [0,a]$ is a geodesic
  curve. Continuing it in $[0,
  \pi]\times [0, a+\pi/2]$ we find a geodesic which hits 
  the boundary $[0,\pi]\times
  \{a+\pi/2\} \cup \{0,\pi\}\times [a,a+\pi/2]$. This lifts
  to a minimal surface $\Sigma^+$ on $\SS^3$ with the metric
  \eqref{e:finalmetric}. Note that the subset of $\Sigma^+$ 
  lying in $S^+$ is a hemisphere. 

  We argue in the same way for $\gamma^-$ and $\gamma_i$. Thus we find
  a sequence of minimal $2$--spheres
  $\Sigma_i$ converging to a singular lamination given by the union of
  $\Sigma^+$, $\Sigma^-$ (lifting of $\gamma^+$ and $\gamma^-$)
  and the strictly stable $2$--sphere
  $\{r=0\}$. Every $\Sigma_i$ contains two hemispheres
  $H^+_i$ and $H^-_i$, lying in $S^+$ and $S^-$.
  All $H^+_i$'s intersect in the great circle given by
  $\{r=a, \phi=\pi/2\}$ (and by symmetry all $H^-_i$'s intersect in
  $\{r=-a, \phi=\pi/2\}$). Thus $\{\Sigma_i\}\cup\{\Sigma_\infty\}$
  can be completed locally to a foliation by great spheres parallel at two
  points.
\end{proof}

\section{Completing the metric and the surfaces; 
proof of Theorem \ref{t:genus}}

In this section we show how to complete the metric (and the 
minimal annuli) constructed near the 
strictly stable $2$-sphere in the previous section.  This will give
Theorem
 \ref{t:genus}, which is 
significantly easier to prove than Theorem \ref{c:gensc}
since we do not require any curvature 
control.  

\begin{figure}[htbp]
\begin{center}
    \begin{minipage}[t]{0.45\linewidth}
    \centerline{\input{pic5.pstex_t}}
    \caption{Metric on the product of an interval with a genus $\bf{g}$ 
 surface with a cylindrical end}
    \label{fig5}
    \end{minipage}
    \hspace{1cm}
    \begin{minipage}[t]{0.45\linewidth}
    \centerline{\input{pic6.pstex_t}}
    \caption{Gluing together the minimal foliation of $N_2$ and the 
minimal foliation of $N_1$}
    \label{fig6}
    \end{minipage}
\end{center}
\end{figure}

\begin{proof}
(Rough sketch of Theorem \ref{t:genus}). 
Let $\Sigma_{\bf{g}}\setminus \{p\}$ be a punctured 
 surface of genus $\bf{g}$ 
equipped with a metric 
which near the puncture $p$ is isometric 
to a flat cylinder.  Let 
$N_1$ be the metric product 
$(\Sigma_{\bf{g}}\setminus U)\times \, ]-\eps,\eps[\,$ 
for some sufficiently small $\eps$; see fig. \ref{fig5}.  
Then $N_1$ is foliated by the minimal surfaces  
$(\Sigma_{\bf{g}}\setminus U)\times \{t\}$.  Let $\Sigma_k$, 
$\Sigma_{\infty}$, $\Sigma_{-\infty}$ be the surfaces constructed in
Proposition \ref{p:neck}. In particular we can assume that they are
lifting of the geodesics $\gamma_k$, $\gamma_\infty$, $\gamma_{-\infty}$ 
of Lemma \ref{l:spirals}.
 
Let $N_2=T_{\nu}(\Sigma_\infty)$ for some sufficiently small $\nu>0$.  
By Lemma \ref{l:spirals} we can assume that part of $N_2$
has a smooth minimal foliation 
of the form $\big\{\SS^1\times\, ]-\delta, \delta [\,
\times\,\{t\}\big\}_{t\in ]-\eps, \eps [}$
where $\SS^1\times\, ]-\delta, \delta [\,  \times \{0\}\subset
\Sigma_\infty$ and, for a sequence $\sigma_k$, 
$\SS^1\times \, ]-\delta, \delta[\,\times \{\sigma_k\} \subset \Sigma_k$. 
The idea is now to glue 
$N_1$ together with $N_2$ along these two foliations while keeping the leaves 
minimal; see fig. \ref{fig6}.  (In Lemma \ref{l:glue} below we will show how to do 
the gluing.) On the other side of $\SS^2\times \{0\}$ we complete the metric 
in the same way except for this time letting the punctured  surface 
have genus $0$.  This gives the desired embedded 
minimal surfaces and the limit lamination in a manifold with boundary 
which is topologically $\Sigma_{\bf{g}}\times \, ]0,1[\,$.  Since 
$\Sigma_{\bf{g}}\times \, ]0,1[\,$ can be topologically embedded into 
$\RR^3$ it  
is now easy to see that the metric can be completed to a metric on the given 
$M$ with the desired property. \end{proof}

To make the construction outlined above precise 
we will need the following two lemmas:

\begin{Lem}  \label{l:neigh}
Let $f$ be a smooth function on $M^3$ with $0$ as a regular value and
let $\Sigma_r=\{f=r\}$ be the level sets of $f$.  
In a tubular neighborhood of $\Sigma_0$ the metric can be written as 
\begin{equation}  \label{e:neigh}
g=k^2(r,\theta)\,dr^2+h(r,\theta)
\end{equation}
where $f(r,\theta)=r$ and $h(r,\cdot)$ is the metric on $\Sigma_r$. 
Moreover, the level sets of $f$
are minimal if and only if $\partial_r \det\, (h)= 0$.
\end{Lem}

\begin{proof} The surface $\Sigma_r$ is minimal if and only if 
  $\text{div}_{\Sigma_r}(\nabla f)=0$.
An easy computation shows that $2 \det\, (h)\, \text{div}_{\Sigma_r}(\nabla
 f)=k\,\partial_r\det\, (h)$ and hence gives the claim.  
\end{proof}

The next lemma shows that we can deform any metric on a product with a
minimal foliation into
the product metric with the product foliation,
while keeping the leaves minimal; see fig. \ref{fig6}.

\begin{Lem}  \label{l:glue}
Let $g$ be a smooth metric of the form \eqr{e:neigh} on 
$\SS^1 \times \, ]0,1+\eps[\, \times \, ]0,1[\,$ for some $\eps>0$ 
and assume that every slice $\SS^1\times \, ]0,1+\eps[\,\times \{t\}$ 
is minimal. Then there is a smooth metric $\tilde{g}$ on 
$\SS^1 \times \, ]0,3[\,\times \, ]0,1[\,$ coinciding with $g$ on
$\SS^1\times \, ]0,1[\,\times \, ]0,1[\,$ and with the product metric on
$\SS^1\times \, ]2,3[\,\times \, ]0,1[\,$ and such that every slice $\SS^1\times
\, ]0,3[\,\times \{t\}$ is minimal.
\end{Lem}

\begin{proof} 
By Lemma \ref{l:neigh} it is sufficient to find a
smooth positive function $\tilde{k}$ and a smooth family of $2$-dimensional 
metrics $\tilde{h}(r,\cdot)$ (both functions of $(r,x,\theta)\in \SS^1\times
\, ]0,3[\,\times \, ]0,1[\,$) with $\partial_r \det\, (\tilde{h})=0$ and
\begin{itemize}
\item $\tilde{k}=k$, $\tilde{h}=h$ on $\SS^1\times \, ]0,1[\,\times \, ]0,1[\,$;
\item $\tilde k=\tilde{h}_{\theta\theta}=\tilde{h}_{xx}=1$,
  $\tilde{h}_{x\theta}=0$ on $\SS^1\times \, ]2,3[\,\times \, ]0,1[\,$\,.
\end{itemize}
(Here and in what follows $h_{\theta\theta}$, $h_{xx}$, 
and $h_{x\theta}=h_{\theta x}$ denote the components of
the metric tensor $h$ in the coordinates $(\theta, x)$; the same
convention is adopted for any other tensor.)
The requirements on $\tilde{k}$ are trivial to satisfy; 
so we only need to construct $\tilde{h}$.
To do that let $\eta: \, ]0,3[\,\to \, ]0, 5/4[\,$ be a smooth function with
\begin{equation}
\eta (x)= \left\{
\begin{array}{lll}
x  & \text{ for }x\in \, ]0,1/2[\, ,\\
3/4 &\text{ for } x\in \, ]1,3[\, ,\\
\end{array}
\right.
\end{equation}
and set
\begin{equation}
h^{(1)}(\theta, x, r)= h (\theta, \eta (x), r) \quad \text{ for }
(\theta,x, r) \in \SS^1\times \, ]0,3[\,\times \, ]0,1[\, .
\end{equation} 
Since $\det\, (h^{(1)}(\theta, x,
r))= \det\, (h(\theta, \eta (x), r))$ 
clearly $\det\, (h^{(1)}(\theta,x,r))$ is
constant in $r$. Next choose a smooth function 
$\varphi: \, ]0,3[\,\to [0,1]$ with  
\begin{equation}
\varphi= \left\{
\begin{array}{lll}
1   & \text{ on }\, ]0,1[\, ,\\
0&\text{ on } \, ]3/2,3[\, .\\
\end{array}
\right.
\end{equation}
Set $h^{(2)}_{\theta\theta}= h^{(1)}_{\theta\theta}$ and
\begin{equation}
h^{(2)}_{x\theta}
(\theta, x, r)= \varphi (x)\, h^{(1)}_{x\theta} (\theta, x,r)\, , \qquad
h^{(2)}_{xx}= h^{(1)}_{xx}\,+\, 
\frac{(1-\varphi^2 (x))\left[h^{(1)}_{x\theta}
(\theta, x,r)\right]^2}{h^{(1)}_{\theta\theta} (\theta,x,r)}\, .
\end{equation}
One easily checks that $h^{(2)}(\cdot,r)$ is a metric for all
$r$ and that $\det (h^{(2)})$ coincides everywhere with
$\det (h^{(1)})$. Thus also $\det (h^{(2)})$ is constant in $r$.
Note that for $x\in \, ]3/2, 3[\,$ the metric $h^{(2)}$ is of the form
\begin{equation}\label{diagonal}
\left(\begin{array}{ll}
h^{(2)}_{\theta\theta} (\theta,x, r) & 0\\
0 & h^{(2)}_{xx} (\theta,x, r)
\end{array}\right)\, .
\end{equation}
Now let $\Phi: \, ]0,3[\,\times \RR^+\times \RR^+\to \RR^+$ be a 
smooth function with 
\begin{equation}
\Phi (x,u,v)= \left\{
\begin{array}{lll}
u  & \text{ for }x\in \, ]0, 3/2[\, ,\\
(uv)^{-1/2}&\text{ for } x\in\, ]2,3[\, .\\
\end{array}
\right.
\end{equation}
Set 
$h^{(3)}_{x\theta}=h^{(2)}_{x\theta}$ and
\begin{eqnarray}
h^{(3)}_{\theta\theta} (\theta,x,r) &=&
\Phi \left(x,\, h^{(2)}_{\theta\theta} (\theta,x,r),\,
h^{(2)}_{xx} (r,x,\theta)\right)\, ,\\
h^{(3)}_{xx}(\theta ,x, r)&=&\frac{h^{(2)}_{\theta\theta}
  (\theta,x,r)\, h^{(2)}_{xx} (r,x,\theta)}
{h^{(3)}_{\theta\theta} (\theta,x,r)}\, .
\end{eqnarray}
Since $\Phi$ takes values in $\RR^+$, $h^{(3)}$ is a well defined
smooth metric. Moreover, we have the identity
$h^{(3)}_{\theta\theta}h^{(3)}_{xx}=h^{(2)}_{\theta\theta} 
h^{(2)}_{xx}$ everywhere.
Since $h^{(3)}$ coincides with
$h^{(2)}$ for $x\in \, ]0,3/2]$
and $h^{(2)}$ is of the form
(\ref{diagonal}) for $x\in \, ]3/2, 3[\,$, 
this yields that $\det (h^{(3)})=\det (h^{(2)})$
everywhere. Note that for 
$x\in \, ]2,3[\,$ we have
$h^{(3)}_{\theta\theta}=h^{(3)}_{xx}$. Moreover, $\partial_r \det
(h^{(3)})=0$ and hence $\partial_r
h^{(3)}_{\theta\theta}(x,\theta, r) =0$ for $x\in\, ]2,3[\,$.
Thus $h^{(3)}$ is of the form
\begin{equation}
\left(
\begin{array}{ll}
\overline{h} (\theta, x) & 0\\
0 & \overline{h} (\theta, x)
\end{array}\right)\, .
\end{equation}
Clearly we can modify $\overline{h}$ for $x\in \, ]5/2, 3[\,$
keeping it as above
for $x\in \, ]2, 5/2]$, positive and
smooth on the whole $\, ]2,3[\,$ and forcing it to be 
identically $1$ in a neighborhood
of $x=3$. This yields the desired metric.
\end{proof}

\begin{proof}
(of Theorem \ref{t:genus}). 
Using Lemma \ref{l:glue} we can now easily carry out the 
gluing outlined in the rough sketch of Theorem \ref{t:genus} above.
\end{proof}

\section{Connected sum construction; proof of Theorem \ref{t:gensc}}

We prove Theorem \ref{t:gensc} by using a connected sum construction.  
When $M$ carries a metric with positive scalar curvature 
this gives a metric 
on $M$ with positive scalar curvature and the desired degenerating sequence 
of minimal surfaces.  For general metrics on general $M$ this
gives a different proof of Theorem \ref{t:genus}.

The connected sum  
is done using in part arguments of \cite{GrLa} and
\cite{ScYa}.  We use the low--tech argument of Gromov and Lawson to
construct an explicit neck connecting two domains in a round
$3$-sphere. (This explicit construction is used when we glue together
minimal surfaces.)  We also use a more high--tech argument of Schoen and
Yau to show that such a metric exists on any $3$-manifold which carries a
metric with positive scalar curvature.  (The result of Schoen and Yau
that we use says that if a $3$-manifold carries a metric of positive scalar
curvature, then the punctured manifold (punctured at a point) has
a metric with $\Scal>0$ and a cylindrical end.) 

Consider again a warped product metric on $\SS^2\times \RR$ of the
form \eqr{e:metric} where $\lambda=\lambda(r)$ is given by
\begin{equation}
\lambda (r)= \left\{
\begin{array}{lll}
-\sin r & \text{ for }r\in [-\pi,-\eps[\, ,\\
\lambda_{GL}(r)&\text{ for } r\in [-\eps,\eps]\, ,\\
\sin r&\text{ for } r\in \,\, ]\eps,\pi]\, .\\
\end{array}
\right.
\end{equation}
Note that the resulting metric is a metric on the $3$-sphere that is
metrically the connected sum of two round unit metrics on the
$3$-sphere by a neck given by the function $\lambda_{GL}$.

\begin{figure}[htbp]
\begin{center}
\input{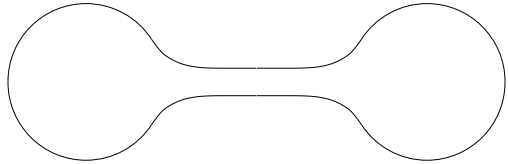}
\caption{The connected sum of two round $\SS^3$'s. The resulting metric has 
positive scalar curvature}
    \label{fig7}
\end{center}
\end{figure}

By section 5 of 
\cite{GrLa} (see also \cite{ScYa}) $\lambda_{GL}$  
can be chosen so that the connected sum still has positive scalar 
curvature for all $\eps>0$; see fig. \ref{fig7}.  
(For completeness we show in Appendix \ref{s:a} 
how to choose
$\lambda_{GL}$ so that the scalar curvature of the warped product is
positive.)  
Call $x$ and $y$ the two points in 
the two copies of $\SS^3$ about where we do the 
connected sum. 

Suppose next that we have two one parameter 
families of minimal surfaces (one in each copy of $\SS^3$).  Suppose that 
one of these families goes through $x$ and the other goes through $y$ and so 
that near $x$, respectively, $y$ the families of minimal surfaces 
are foliations by great $2$-spheres.  
We show in Lemma \ref{l:scalglue} below that when we take 
the connected sum of the two $\SS^3$'s by a neck as above, 
then we can glue the minimal surfaces 
in one of the two $3$-spheres together with the minimal surfaces in the other 
$3$-sphere keeping the surfaces minimal through the neck. In 
Lemma \ref{l:model} 
below 
we then show that we can find a metric on $\SS^3$ with positive 
scalar curvature and with a family of embedded minimal tori going through a 
point as a foliation by great spheres on a round unit $\SS^3$.  Taking 
the connected sum of ${\bf{g}}$ copies of 
this metric on $\SS^3$ with the metric on $\SS^3$ and minimal spheres 
constructed in Proposition \ref{p:basiclam} will then  
prove Theorem  \ref{t:gensc} when $M=\SS^3$.    
Finally, taking the connected sum (using theorem 4 of \cite{ScYa}) 
with a general $M^3$ we get 
Theorem \ref{t:gensc}.

\begin{figure}[htbp]
\begin{center}
\input{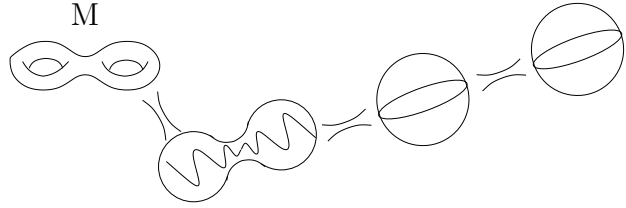}
\caption{Connected sum of one--parameter families of tori in $\bf{g}$ $3$--spheres
with the desired degeneration and $M$}
    \label{fig8}
\end{center}
\end{figure}

In Lemma \ref{l:scalglue} we are able to glue only 
foliations which are parallel (see Definition \ref{d:parallel}). 

\begin{Lem}  \label{l:scalglue}
Let $\Omega_1$, $\Omega_2$ be two open subsets of the round unit $\SS^3$
with $x_i\in \Omega_i$.  Suppose that $\cF_1$ and $\cF_2$ are
foliations by great spheres parallel at $x_1$ and $x_2$ respectively.
In $\Omega_1\#\Omega_2$ we can connect the central leaves and the
ones nearby keeping them minimal and $\Scal>0$.
\end{Lem}

\begin{figure}[htbp]
\begin{center}
    \begin{minipage}[t]{0.45\linewidth}
    \input{pic9.pstex_t}
    \caption{The connected sum construction}
    \label{fig9}
    \end{minipage}
    \begin{minipage}[t]{0.45\linewidth}
    \input{pic10.pstex_t}
    \caption{The rectangle $Ap$ (corresponding to zone A). The metric 
\eqref{e:modifiedplan} is the standard metric on $\SS^2$}
    \label{fig10}
    \end{minipage}
\end{center}
\end{figure}

\begin{proof} Let $(\phi, \theta, r)$ be spherical 
coordinates on $\SS^3$ centered at $x_1\in \Omega_1$.
The standard metric is $dr^2 +\sin^2 r\, (d\phi^2 + \sin^2 \phi\,
d\theta^2)$. Endow the square $[0, \pi]\times [0, \pi]$
with the degenerate metric $\sin^2 r \sin^2 \phi\, ( dr^2 + \sin^2 r\,
d\phi^2)$. Clearly the
geodesics passing through $(\pi/2, \pi/2)$ lift to great spheres
parallel at $x_1$, with $\{\phi=\pi/2\}$ the central leaf 
containing $x_1$. We do the same at $x_2\in \Omega_2$.
  
  The construction outlined in Appendix A shows
  that we can replace the balls $B_{\eps} (x_1)\subset \Omega_1$ 
  and $B_{\eps}
  (x_2)\subset \Omega_2$ with two hyperbolic necks and then connect the
  two necks with a cylinder $\SS^2\times \, ]-K,K[\,$. 
  More precisely, the construction gives a
  metric on $\SS^2\times \, ]-K_1-\eps, K_1+\eps[\,$ of the form 
  \begin{equation}\label{e:warped}
  dr^2 + \lambda^2 
  (r)\, (d\phi^2 + \sin^2 \phi\, d\theta^2)
  \end{equation}
   where, see fig. \ref{fig9},
  \begin{itemize}
  \item[-] $\lambda (r)=\lambda (-r)$ and $\lambda'(r)\geq 0$ for $r\geq 0$;
  \item[-] $\lambda (r)=\sin (r- (K_1-\eps))$ for $r\in [K_1, K_1+\eps]$;
  \item[-] $\lambda$ has a hyperbolic behavior on $[K, K_1]$ and is
    constant $R$ on $[0, K]$.
  \end{itemize}
  \noindent Note that we can make the cylindrical tube as long as we
  want (in particular we can assume that $K> \pi R/2$).
  The
  coordinates $(\phi, \theta)$ have been chosen
  in such a way that the leaves of the minimal 
  foliations 
  are lifting of two families of geodesic segments in $[0, \pi]\times
  \, ]-K_1-\eps, K_1+\eps[\,$ with the corresponding degenerate metric
  \eqref{e:cube}. We can continue our geodesic
  segments throughout the whole strip $[-\pi, \pi]\times
  \, ]-K_1-\eps, K_1+\eps[\,$. They do not hit the boundary
  lines $\{\phi=0\}$, $\{\phi= \pi\}$ and they give two one--parameter families of
  geodesics $\cG_1$ and $\cG_2$, which lift to mininal
  surfaces in the metric \eqref{e:warped}. 
  These minimal surfaces are all cylinders: Their boundaries
  are two circles lying on $\SS^2\times \{K_1+\eps\}$ and $\SS^2\times
  \{-K_1-\eps\}$.

  Note the following:
  \begin{itemize}
  \item[(i)] The two central leaves
    are lifting of the ``central'' geodesic $\gamma_c=\{\phi=\pi/2\}$, 
    hence they naturally connect.
  \item[(ii)] Both $\cG_i$ are symmetric around $\gamma_c$, 
    i.e., if $\gamma = \{(\phi (t), r(t))\}_{t\in ]a,b[}$ 
   lies in $\cG_i$, then so does $\{(\pi-\phi (t), r(t))\}_{t\in ]a,b[}$. 
  \item[(iii)] If $\{(\phi (t), r(t))\}_{t\in ]a,b[}$ lies in $\cG_1$, then
    $\{(\phi (t), -r(t))\}_{t\in ]a,b[}$ lies in $\cG_2$.
  \end{itemize}
  Together (ii) and (iii) give
\begin{equation}\label{e:symmetry}
\mbox{if }\{(\phi(t), r(t))\}_{t\in ]a,b[} 
\text{ lies in }\cG_1,\text{ then }
\{(\pi-\phi(t),-r(t))\}_{t\in ]a,b[}
\text{ lies in }\cG_2\, .
\end{equation}

  For $\eps>0$ sufficiently small we can modify the metric in
  $\SS^2\times [-K,K]$ so that:
  \begin{itemize}
  \item[(a)] It has $\Scal>0$ and is of the form 
    \begin{equation}\label{e:specform}
    k^2 (r, \phi)\,dr^2 + R^2\,
    (d\phi^2 + g^2 (r, \phi)\, d\theta^2)\, .
    \end{equation}
  \item[(b)] $k(r, \phi)= k(-r, \phi)= k(r, \pi-\phi)$ and the same is
    true for $g$.
  \item[(c)] In an $\eps$--neighborhood of $\{\phi=\pi/2, r\in
    [-\pi R/2, \pi R/2]\}$ (``zone A'' in fig. \ref{fig9})
the metric is 
   \begin{equation}\label{e:rsphere}
    \sin^2 \phi\,
    dr^2 + R^2\, (d\phi^2 + d\theta^2)\, .
    \end{equation} 
  \item[(d)] The cylinder $\{\phi=\pi/2\}$ remains a minimal surface.
  \end{itemize}
  That this modification is possible can be shown in the same
  way as Lemma \ref{l:model} (cf. the
  second step of the proof). We give the details of this at the end. 
  We first show how in the new metric the two foliations connect
  nearby $\{\phi=\pi/2\}$. 

  By (a) we can apply the discussion of Section
  1. The families
  $\cG_1$ and $\cG_2$ become two new families of curves
  $\cG'_1$ and $\cG'_2$, which are geodesics in the modified
  metric 
  \begin{equation}\label{e:modifiedplan}
  R^2\, g^2 (\phi, r)\, (k^2 (\phi,r)\, dr^2 + R^2\, d\phi^2)\, .
  \end{equation}
  $\cG'_1$ coincides with
  $\cG_1$ for $r> K$, whereas $\cG'_2$ coincides with
  $\cG_2$ for $r<-K$. Moreover, $\gamma_c$ is a geodesic also for
  \eqref{e:modifiedplan} and lies in both $\cG'_1$ and $\cG'_2$. 
  By continuity of the dependence on the initial data, all the curves
  of $\cG'_1$ which in
  $\{r>K_1\}$ start sufficiently near $\gamma_c$ 
  never leave its $\eps$--neighborhood (and are all graphs of
  functions of $r$).
  
  In the rectangle 
  $Ap=\, ]\pi/2-\eps, \pi/2+\eps[\,\times \, ]-\pi R/2, \pi R/2[\, $
  the metric \eqref{e:modifiedplan} is given by
  \begin{equation}\label{e:isthesphere}
  R^2 \sin^2 \phi\, dr^2 + R^4\, d\phi^2\, .
  \end{equation}
  Note that \eqref{e:isthesphere}
  is the metric on the round $2$--sphere of radius
  $R^2$ and $\gamma_c\cap Ap$ is half of a great circle. 
  
  Now take a $\gamma\in \cG'_1$
  which intersects $\{r=\pi R/2\}$ 
  transversally and leaves $Ap$ crossing $\{r=-\pi R/2\}$. 
  Also $\gamma\cap Ap$ is half of a great circle. 
  It is easy to check that the crossings of $\gamma$ with $\{r=\pi R/2\}$ 
  and $\{r=-\pi R/2\}$ are two antipodal points.
  Thus if $\gamma$ crosses $\{r=\pi R/2\}$ at 
  $\phi=\phi_0$ with angle $\delta$, it 
  crosses $\{r=-\pi R/2\}$ at $\phi=\pi-\phi_0$ with
  angle $-\delta$ (see fig. \ref{fig10}). By (b), the families $\cG'_1$ and
  $\cG'_2$ satisfy condition \eqref{e:symmetry}. Thus there is a
  geodesic in $\cG'_2$ which crosses $\{r=-\pi R/2\}$ at
  $\phi=\pi-\phi_0$ with angle $-\delta$. This geodesic connects with
  $\gamma$.

\medskip

\noindent
\underline{The modified metric}:  
We complete the proof by showing how to construct the modified
metric. Straightforward computations give that
for a metric of the form \eqref{e:specform}
\begin{equation}\label{e:scalmod}
\Scal= -\left(\, \frac{k_{\phi\phi}}{kR^2} + \frac{g_{\phi\phi}}{g
    R^2}\,\right)\, -\,
\left(\,\frac{g_{rr}}{gk^2}
+ \frac{g_\phi k_\phi} {g k R^2} - \frac{g_r k_r}{g k^3}\, \right)\, .
\end{equation}
Fix a bump function $\varphi: [0, K]\to [0,1]$ which is $0$
in a neighborhood of $K$ and is 1 in a neighborhood of $[0, \pi
R/2]$. Let $C$ be a constant such that $|\varphi'|, |\varphi''|\leq
C$. 

It is easy to check that for any $\eps>0$ we can find functions
$\tilde{g}, \tilde{k}:[\pi/2,\pi]\to [0,1]$ such that
\begin{itemize}
\item[($\alpha$)] $\tilde{k} (\phi)= \sin \phi$ in a neighborhood $I$ of
  $\pi/2$ and is 1 outside another neighborhood.
\item[($\beta$)] $\tilde{g} (\phi)= \sin \phi$ {\em outside} $I$ and
  is $1$ in a smaller neighborhood of $\pi/2$.
\item[($\gamma$)] $|\tilde{g}-1|\leq \eps$ where $\tilde{g}$ differs
  from sine; $|\tilde{k}-1|, |\tilde{k}'|,
  |\tilde{g}'|\leq \eps$ and $\tilde{k}'', \tilde{g}''\leq \eps$ 
everywhere.  
\end{itemize}
The functions $k$ and $g$ are then given by 
\begin{equation}
g(\phi, r)= \varphi (r) \,\tilde{g}(\phi) + (1-\varphi
(r))\,\sin \phi\, , \qquad
k(\phi, r)= \varphi (r) \,\tilde{k}(\phi) + (1-\varphi(r))
\end{equation}
on $[\pi/2, \pi]\times [0, K]$ and we extend them by symmetry to $[0,
\pi]\times [-K, K]$. The resulting metric is smooth and coincides with
the product outside a neighborhood of
$\{\phi=\pi/2, r\in [-\pi R/2, \pi R/2]\}$
in $\SS^2\times [-K, K]$. Clearly, $k$ and $g$ satisfy (b),
(c), and, by Lemma \ref{e:neigh}, (d).

To complete the proof we need to show that the scalar curvature is
positive where the metric differs from the standard product. 
It is easy to check that $|\partial_r k|, \partial_\phi
k|, |\partial_ r g|, |\partial_\phi g|, |\partial_{rr} k| \leq
C\eps$. Thus, for $\eps$ small,
\begin{equation}\label{e:firstest}
\left|\,
\frac{g_{rr}}{gk^2}
+ \frac{g_\phi
  k_\phi} {g k R^2} - \frac{g_r k_r}{g k^3}\, \right|\leq 2C\eps
(1+R^{-4}+ R^{-2})\, .
\end{equation}
Moreover, if $\eps$ is sufficiently small, 
($\alpha$), ($\beta$), and the inequalities $\tilde{k}'',
\tilde{g}''\leq \eps$ give that 
\begin{equation}
- \frac{g_{\phi\phi}}{g R^2}
- \frac{k_{\phi\phi}}{kR^2}\geq \frac{1}{2R^2}\, .
\end{equation}
Since $\eps$ can be chosen arbitrarily this completes the proof.
\end{proof}

On $\SS^2\times \SS^1$ with the product metric, the great circles on
$\SS^2$ times $\SS^1$   
give a one parameter family of minimal (intrinsically flat) tori. The 
next lemma shows that we can deform this example into a one parameter 
family of embedded minimal tori on $\SS^3$ with a metric with positive 
scalar curvature and so that in a neighborhood of some point the metric 
has constant sectional curvature $1$ and the tori pass through as 
parallel great $2$-spheres. 
(The proof of this lemma is postponed to Appendix \ref{s:b}.)

\begin{Lem}\label{l:model} 
On $\SS^3$ there exists a metric with 
$\Scal>0$
and a family of minimal tori 
$\{T_\delta\}_{\delta\in \, ]-1,1[}$ such 
that in a neighborhood of two antipodal points $x$ and $y$
the metric coincides with the 
round unit metric and $\{T_\delta\}$ with a foliation by great
$2$--spheres parallel at $x$ and $y$.
\end{Lem}

\begin{proof}
(of Theorem \ref{t:gensc}). 
Let $M_{\text{tor}}$ be the metric on
$\SS^3$ given by Lemma \ref{l:model} and let $M_{\text{sing}}$ be the
metric on $\SS^3$ given by Proposition \ref{p:basiclam}.  
By Lemma \ref{l:scalglue}
$\#_{i=1}^{\bf{g}}M_{\text{tor}}\#M_{\text{sing}}$ gives a metric on
$\SS^3$ and a sequence of embedded minimal surfaces of genus $\bf{g}$
with the desired properties. (Here all necks 
are attached at points where the sectional curvatures are constant.)
By theorem 4 of Schoen-Yau, \cite{ScYa}, 
there exists a metric on $M^3$ with positive scalar curvature
and a cylindrical end. Connecting this metric
with the metric on the $3$-sphere constructed above completes the proof.
(The last neck connects the cylindrical end with an open set of
the $3$-sphere where the sectional curvatures are constant.)
\end{proof}

\section{Metrics on $\SS^2\times \SS^1$}

\begin{proof}
(of Theorem \ref{t:product}).  This is essentially proven in \cite{HaNoRu} 
although not recorded there.   Namely, similarly to Proposition \ref{p:neck} 
consider the 
degenerate metric 
\begin{equation}  \label{e:stproduct}
ds^2=\sin^2\phi\,(dr^2+d\phi^2)
\end{equation}
on the cylinder $[0,\pi]\times \SS^1$.  Geodesics in 
this metric lift to minimal surfaces on the product $\SS^2\times
\SS^1$ and simple closed geodesics lift to embedded minimal tori.  
By lemma 2.1 
of \cite{HaNoRu} geodesics  
in \eqr{e:stproduct} are periodic (in $r$) 
and as the angle that they make with the 
geodesic $\{r=0\}$ goes to zero the period in $r$ goes to zero.  
Moreover, it follows easily that the 
period is continuous as a function of the angle.  Combining these 
facts is easily seen to give 
that there are simple closed 
geodesics on 
the cylinder with arbitrarily small period in $r$ and that these 
converge to the foliation of the cylinder by the parallel geodesics 
$\{r=\mbox{ constant}\}$.  
Lifting these simple closed geodesics to $\SS^2\times \SS^1$ 
gives the desired sequence of 
embedded minimal tori.
\end{proof}

\begin{Rem} Arbitrary close to the product metric on 
$\SS^2\times \SS^1$ we can also find a metric which has 
  a sequence of embedded minimal tori converging to a singular lamination
  of the type of Theorem \ref{t:genus}. Indeed we choose on
  $\SS^2\times \RR$ a metric of the form (\ref{e:metric}) where
  $\lambda (r)$ is symmetric, equal to
  $\cosh (\eps r)$ for $r\in \, ]-1,1[\,$ and
  constant on $\, ]-\infty, 2]\cup [2, \infty[\,$. Consider 
  on the strip $[0,\pi]\times \RR$ the degenerate metric
  whose geodesics lift to
  minimal surfaces on $\SS^2\times \RR$.
  By Lemma \ref{l:spirals} for any given $r_0>2$ there is a sequence of
  geodesics $\gamma_i$ which all pass through $(\pi/2, r_0)$, $(\pi/2, 0)$ 
  and $(\pi/2, -r_0)$ and which converges to
  a lamination consisting of $\{r=0\}$ and two infinite geodesics
  $\gamma_\infty$ and $\gamma_{-\infty}$ spiraling into it. We now identify
  the lines $\{r=r_0\}$ with $\{r=-r_0\}$ on the strip and the spheres
  $\{r=r_0\}$ and $\{r=-r_0\}$ in $\SS^2\times \RR$. Thus we obtain a
  smooth metric on $\SS^2\times \SS^1$ and a degenerate metric on
  $[0,\pi]\times \SS^1$ whose geodesics lift to minimal surfaces in
  $\SS^2\times \SS^1$. Because of the symmetry of our
  construction the geodesics $\gamma_i$ generate simple closed
  geodesics in $[0, \pi] \times \SS^1$ and $\gamma_\infty$ and
  $\gamma_{-\infty}$ smoothly glue themselves forming an infinite
  geodesic spiraling into $\{r=0\}$ from both sides. These
  geodesics lift to the desired minimal surfaces in $\SS^2\times
  \SS^1$.
\end{Rem}


\section{More than one strictly stable $2$-sphere with singularities}  
\label{s:metrics}

The proof of Theorem \ref{t:genus} easily generalizes to show that for any
given integer $n>0$ we can find a limit lamination which is 
singular at $n$ pairs
of points, where the pairs of points lie on $n$ disjoint strictly stable 
$2$-spheres.  That is:

\begin{Thm}   \label{t:genusa}
On any $3$-manifold, $M^3$, and for any nonnegative integer $\bf{g}$,
and any positive integer $n$ 
there exists a metric on $M$ and a sequence of embedded minimal 
surfaces of 
genus $\bf{g}$ with Morse index going to infinity and which converges to a 
singular (minimal) lamination $\cL$.  This can be done so that the singular 
set of $\cL$ consists of pairs of points lying on $n$ leaves 
which are strictly 
stable $2$-spheres.
\end{Thm}

\begin{proof}
Let $M_{\delta}$ be a $\delta$-tubular neighborhood of 
the strictly stable $2$-sphere in
the metric on $\SS^2\times \RR$ given by Proposition \ref{p:neck}.
Using  Lemma
\ref{l:glue} glue $n$ copies of $M_{\delta}$ together
along the minimal leaves of the foliation near the boundary coming
from 
the $\Sigma_i$'s 
and $\Sigma_\infty$ while keeping the leaves minimal.  The desired 
metric can now be obtained by completing this
metric to a metric on $M$ using Lemma \ref{l:glue} as in the proof of
Theorem \ref{t:genus}.
\end{proof}

Likewise we can easily generalize
Theorem \ref{c:gensc} to:

\begin{Thm}   \label{t:gensca}
If $M^3$, $\bf{g}$ are as in Theorem \ref{c:gensc}, and $n$ is a 
positive integer, then 
there is a metric with $\Scal_M>0$ which has a singular 
lamination and a sequence of embedded 
minimal surfaces of genus ${\bf{g}}$ as in Theorem \ref{t:genusa}.
\end{Thm}

\begin{proof} We use the same ideas of the proof of Proposition
  \ref{p:basiclam} to glue $n$ hyperbolic necks and $2$ halves of
  standard $\SS^3$. Thus we produce a metric on $\SS^3$ with positive
  scalar curvature and with a sequence of embedded minimal spheres
  which converge to a singular lamination containing $n$ strictly
  stable $2$--spheres and which pass through two points as parallel
  great spheres. We can use the connected sum construction of
  Section 3 to complete the proof.
\end{proof}

\section{The space of minimal annuli limiting a strictly
  stable $2$-sphere}         \label{s:speculation}

Recall the following theorem from \cite{CH2} (here $T_1M$ is the unit
tangent bundle):

\begin{Thm} \label{c:stablemfldf}
\cite{CH2}.  
Let $M^2$ be an orientable surface and $\gamma\subset M$  
be a simple closed and strictly stable geodesic.
Then there are four ``circles'' of
noncompact
geodesics limiting on $\gamma$.  That is, on each side
of $\gamma$ in $M$,
and for each
orientation of $\gamma$ there is a $C^1$ map
$\SS^{1}\rightarrow T_{1}M$ which gives a bijection
between the circle $\SS^{1}$ and the set of geodesics
$\ell$ with $\overline{\cap_{t>0}\ell |[t,\infty[}=\gamma$ 
which limit on $\gamma$ from the given side of $M$ with the given orientation.
\end{Thm}

Motivated by this theorem one is tempted to ask:

\vskip2mm
\noindent
{\bf{Question 1}}:  Let $M^3$ be an orientable $3$-manifold 
and $\Gamma\subset M$ a strictly stable embedded $2$-sphere.  
Does there exist a map from the space of
noncompact embedded minimal annuli   
in $T$ limiting $\Gamma$ and into $\SS^2\times \SS^1\times \ZZ/2\ZZ$?
\vskip2mm

A particular case of the reverse of this question is:

\begin{Pro}\label{p:param}
Let $M=\SS^2\times \RR$ with a metric \eqr{e:metric} where 
$\lambda'(0)=0$, and $\lambda''(0)>0$.  In a neighborhood of the
strictly stable $2$-sphere $\Gamma=\SS^2\times \{0\}$ there are at
least two
``$2$-spheres'' of $\SS^1$ invariant  noncompact
minimal annuli limiting on $\Gamma$. That is, on each side
of $\Gamma$ in $M$,
there is a continuous map
from $\SS^{2}\times \SS^1$ 
to the set of minimal annuli
$\Sigma$ which are the preimages of geodesics $\sigma$ in
\eqr{e:cube} 
with $\overline{\cap_{t>0}\sigma | [t,\infty[}=\{r=0\}$ 
from the given side.
\end{Pro}

\begin{proof}
For each $x\in \SS^2$ we can use spherical coordinates $(\phi,\theta)$ 
centered at $x$ and consider 
the corresponding degenerate metric \eqr{e:cube}  
on the strip. By Lemma \ref{l:spirals} for
every $\rho>0$ we can find a geodesic passing through $(\pi/2, \rho)$
with a given fixed orientation and
spiraling into $\{r=0\}$. This gives a circle worth of annuli.
Varying $x$ gives the claim.
\end{proof}
A weaker question is:

\vskip2mm
\noindent
{\bf{Question 2}}:  Let $M^3$, $\Gamma$ be as in Question 1.   What
is the tangent space of the set of all noncompact embedded minimal annuli
limiting $\Gamma$?  In particular, for each such minimal annulus, what
can be said about the dimension of the space of Jacobi fields that
come from a variation of such annuli?
\vskip2mm

\vskip2mm
In view of \cite{CH2} and \cite{CM3} it seems plausible that if the
answer to Question 1 is yes, then the
following should be the case:

\vskip2mm
\noindent
{\bf{Question 3}}:  Let $M^3$ be an orientable $3$-manifold 
with a generic metric with positive scalar curvature.   Is every 
singular minimal lamination, which is
the limit of a sequence of embedded minimal surfaces of a given fixed
genus, singular along at most one strictly stable $2$-sphere (which
is a leaf of the lamination)? 
\vskip2mm

\appendix

\section{The connecting neck}   \label{s:a}

For completeness we show here how to connect sum two
round unit $\SS^3$'s by a thin neck so that the resulting metric has 
positive scalar curvature
everywhere (cf. section 5 of \cite{GrLa}).  
Fix one of the spheres and a point $x$ on it. Choose spherical
coordinates centered at $x$. In these coordinates the metric is
given by $dr^2 + \sin^2 r\, (d\phi^2 + \sin^2 \phi\, d\theta^2)$.
Starting from $r=\eps$ we will replace $\sin r$ with a function
$\lambda$ and modify the metric as $dr^2 + \lambda^2 (r)\,(d\phi^2 +
\sin^2 \phi\, d\theta^2)$. Shift the coordinate $r$ so that the
replacement of $\sin$ (and hence the neck) starts
at $\{r=0\}$: thus our function $\lambda$ is given
by $\sin (\eps +r)$ on $\{r>0\}$. Hence, 
$\lambda (0)= \sin \eps$ and
$\lambda' (0)=\cos \eps$. Our goal is to continue $\lambda$ in $C^{1,1}$,
while keeping $\Scal>0$ and reaching $\lambda' (-K)=0$ for some
$K>0$ (keeping $\lambda$ positive, so the metric
is not degenerate). For $r<-K$ let $\lambda$ be constant:
hence our metric turns out to be the product of a half--line with the 
round $2$--sphere of radius $\lambda (-K)$.

We make the same construction for the other unit sphere and
then glue the two cylindrical parts. This gives a
$C^{1,1}$ metric which
can be smoothed in a standard way to a
metric with $\Scal>0$.

We will construct $\lambda$ and $-K$ so to have $\lambda''\geq 0$
on $[-K,0]$. Thus 
\begin{equation}\label{e:derivative}
0< \lambda'(r)\leq \lambda' (0)=\cos \eps
\quad \mbox{ for $r\in \, ]-K,0]$\,
  ,}
\end{equation}
\begin{equation}\label{e:invert}
\mbox{$\lambda$ is invertible on 
$\, ]-K,0]$ with $\lambda^{-1}=\alpha$\, .}
\end{equation} 
By \eqref{e:derivative}, 
 $(\lambda')^2\leq (1-\eta)$ on $[-K,0]$ for some $\eta>0$.  
Constructing $\lambda$ in 
this way we will have by \eqr{e:scalwrap}
\begin{equation}
\Scal_M\geq -2 \frac{\lambda''}{\lambda}+\frac{\eta}{\lambda^2}\, .
\end{equation}
Thus we need to find $\lambda$ satisfying
\begin{equation}
\frac{\eta}{4\lambda}\geq \lambda''\geq 0 
\quad \textrm{ on $[-K,0]$}\, ,
\quad\lambda' (-K)=0
\, , \quad \lambda'(0)=\cos \eps\, ,
\quad
\lambda (0)=\sin \eps\,  .
\end{equation}
To do this we solve backward in time the ODE $\lambda''= \eta/
(4\lambda)$ and prove that there is $K>0$ large 
enough so that $\lambda' (-K)=0$ somewhere and $\lambda$
is positive on $[-K, 0]$.  Indeed set $-K= \inf \{t|
\lambda (t)> 0\text{ and }\lambda'(t)>0\}$. 
We claim that if $-K>-\infty$, then
$\lambda (-K)>0$.
If not, then we get the contradiction
\begin{eqnarray}
\lambda' (0)&\geq& \lambda' (0)-\lambda' (-K)=\int_{-K}^0 \lambda'' (t)
\, dt
=\int_{-K}^0 \frac{dt}{4\lambda (t)}\\
&\stackrel{\eqref{e:invert}}{=}&
\frac{1}{4}\int^{\lambda (0)}_0 \frac{d\tau}{\tau\,\lambda'(\alpha (\tau))}
\stackrel{\eqref{e:derivative}}{\geq}
\frac{1}{4\cos \eps}\int_0^{\lambda (0)} \frac{d\tau}{\tau} =
\infty\, .\notag
\end{eqnarray}
Thus, either $-K=-\infty$ or it is
finite and $\lambda' (-K)=0$. In the first case we would have 
\begin{equation}
\lambda' (0)\geq \lambda' (0) - \lim_{x \to \infty} \lambda' (x)=
\int_{-\infty}^0 \lambda'' (t) dt
= \int_{-\infty}^0 \frac{dt}{4\lambda (t)}\geq \frac{1}{4\lambda (0)}
\int_{-\infty}^0 dt =\infty\, .
\end{equation}
This gives a
contradiction; thus $-K>-\infty$ and $\lambda' (-K)=0$.

\section{Proof of Lemma \ref{l:model}}   \label{s:b}

\begin{proof} \underline{A metric with $\Scal>0$ on $\SS^3$ containing
totally geodesic tori}.

We first exhibit a metric on $\SS^3$ with positive scalar curvature 
containing a neighborhood of a totally geodesic torus (given by a
great circle times $\SS^1$) in $\SS^2\times \SS^1$ with the product metric.
The induced metric on a tubular neighborhood $T$ of such a totally geodesic
torus is
\begin{equation}\label{e:standard}
dx^2 + \cos^2 x\, dy^2 + dz^2\, ,
\end{equation}
where $(x,z)$ are coordinates on the torus.
Choose two functions $f,k:[-\pi/2,\pi/2]\to
[0,1]$ and $a\in (0, \pi)$ such that
\begin{itemize}
\item[-] both are positive on $\, ]-\pi/2, -\pi/2[\,$;
\item[-] $f$ coincides with cosine on $\, ]-\pi/2, -2a]$ and 
is $1$ on $\, ]-a, \pi/2[\,$;
\item[-] $k$ coincides with cosine on $[-2a, \pi/2]$ and is constant in a
neighborhood of $-\pi/2$;
\item[-] $f''\leq 0$ and $k''/k\leq 1/4$. 
\end{itemize}
All these conditions can be satisfied
provided $a$ is sufficiently small. We now take $M=[-\pi/2, \pi/2]\times
\SS^1\times \SS^1$ with the metric
$dx^2 + g^2 (x)dy^2+ f^2 (x) dz^2$. Note that the scalar curvature of
this metric is $-k''/k-f''/f$. 
Define on $M$ the equivalence relation $(-\pi/2,x,y)\approx
(-\pi/2, x,z)$ and $(\pi/2, y,x)\approx (\pi/2,z,x)$. 
$M/\approx$ is obtained by gluing 
two solid tori along their boundary (exchanging parallels and meridians)
and thus it is a $3$--sphere.
The metric on $M/\approx$ is smooth
and has positive scalar curvature.

\medskip

\underline{Deforming parts of minimal tori into parts of great spheres}.

The standard metric on $\SS^3$ is
\begin{equation}\label{e:Sphmetric}
\cos^2 \phi\, \cos^2 \theta\, dr^2 + d\phi^2 + \cos^2 \phi\,
d\theta^2\quad \big(=
\cos^2 \phi\, (\cos^2 \theta\, dr^2 + d\theta^2)+ d\phi^2\big)\, ,
\end{equation}
where $\{r=\mbox{ constant}\}$ give a one parameter
families of great spheres parallel in $(0,0,0)$ (see Definition
\ref{d:parallel}).

The product metric on $\SS^2\times \SS^1$ is given by
\begin{equation}\label{e:standard2}
\cos^2 \phi \,dr^2 + d\phi^2 + d\theta^2\, .
\end{equation}
(Here $(\phi,r)$ are spherical coordinates on $\SS^2$ and
$\theta$ is the standard coordinate on $\SS^1$.) Note
that the level sets $\{r=\mbox{ constant}\}$ is a one parameter family
of totally geodesic tori. We modify
$\eqref{e:standard2}$ in a neighborhood of $(0,0,0)$ so that in a
smaller neighborhood the metric is \eqref{e:Sphmetric}, 
the scalar curvature is everywhere positive and
$\{r=c\}$ is a minimal torus (for all $c$ sufficiently small). Notice that
we can do the same modification around the point $(0,\pi, 0)$. 

We first take care of the term in front of $d \theta^2$. 
We can find a function $k$ which is 
$1$ outside a neighborhood of $(0,0,0)$,
coincides with $\cos \phi$
in a smaller neighborhood and does not depend on $r$ if $r$ 
is sufficiently
small. Moreover, for every $\eta$ we can find such a $k$ so that: 
\begin{itemize}
\item[(a)] $|k-1|$ and the norm of all first and second
  partial derivatives of $k$ but $\partial_{\phi\phi} k$ are less than
  $\eta$.
\item[(b)] $\partial_{\phi\phi} k \leq \eta$.  
\end{itemize}
Since for $r$ sufficiently small $k$ does not depend on $r$, by Lemma 
\ref{l:neigh} the leaves
$\{r=\mbox{ constant}\}$ are still minimal in the modified metric
(for $r$ small). 
Moreover, the scalar curvature of the metric $\cos^2 \phi\, dr^2 +
d\phi^2 + k^2 (\phi,r,\theta)\, d\theta^2$ is
\begin{equation}
1 - \frac{\partial_{\phi\phi} k}{k}- \frac{\partial_{rr} k}{k\, \cos^2
 \phi} -
 \tan \phi\, \frac{\partial_\phi k}{k}\, .
\end{equation}
Thus $k$ can be chosen so that the scalar curvature remains
positive. In a neighborhood of $(0,0,0)$ our new metric is 
$\cos^2 \phi\,
dr^2 + d\phi^2 + \cos^2 \phi\, d\theta^2$.
Similarly, we can further modify the metric in a smaller neighborhood so to
adjust the term in front of $dr^2$.

We conclude the proof by constructing the function $k$. Take a smooth cut--off
function $\varphi: \, ]-\delta, \delta[\,\to [0,1]$ which is $0$ in a
neighborhood of $-\delta$ and $\delta$, and $1$ in a neighborhood of
$0$. For some constant $C>1$ we will have $|\varphi'|,
|\varphi''|\leq C$. Next choose a function $\tilde{k}: (-\delta,
\delta) \to [0,1]$ equal to 1 in a neighborhood of $-\delta$ and
$\delta$, equal to $\cos$ in a neighborhood of $0$ and such that
$|k-1|, |k'|\leq \eta/C$ and $k''\leq \eta/C$. (This is
possible since $\cos (0)=1$ and $(\cos)' (0)=-\sin (0)=0$.)
Set 
\begin{equation}
k(r, \phi, \theta)=[1-\varphi (\theta)]\,  +  \varphi (\theta)\, 
[(1-\varphi (r))\, + \, \varphi (r) \tilde{k} (\phi)]\, .
\end{equation}
Clearly, 
$k$ is $1$ in a neighborhood of
the boundary of $[-\delta, \delta]^3$ and does not depend on $r$ if $r$ is
  sufficiently small. Moreover, $|k-1|\leq
  \eta /C$ and 
\begin{equation}
\partial_r k (r, \phi, \theta) = \varphi (\theta)\, \varphi' (r)\, (\tilde{k}
(\phi)-1)\, .
\end{equation}
Hence $|\partial_r k|\leq |\varphi'||\tilde{k} -1|\leq \eta$. 
We argue similarly for all first and second partial derivatives except for
$\partial_{\phi\phi} k$. Finally, $\partial_{\phi\phi} k (r, \phi, \theta)
= \varphi (\theta)\, \varphi (r)\, \tilde{k}'' (\phi)\, \leq \eta$.
\end{proof}

\end{document}